\documentstyle[art10]{article}
\title{Projective normality of special scrolls II.}
\author{Luis Fuentes Garc\'{\i}a\thanks{Supported by an F.P.U.
fellowship of Spanish Government}
 \and{Manuel Pedreira P\'erez}
 }
\date{}

\newtheorem{teo}{Theorem}[section]

\newtheorem{prop}[teo]{Proposition}

\newtheorem{lemma}[teo]{Lemma}

\newtheorem{rem}[teo]{Remark}

\font\euf=eufm10 at 12pt

\def\g2{\pi}
\def\L{{\cal L}}

\def\e{\mbox{\euf e}}

\def\b{\mbox{\euf b}}

\def\a{\mbox{\euf a}}

\def\k{{\cal K}}
\def\K{{\cal K}}

\def\P{{\bf P}}

\newcommand\E{{\cal E}}
\newcommand\Te{{\cal O}}

\def\qed{\hspace{\fill}$\rule{2mm}{2mm}$}
\newcommand\lrw{\longrightarrow}

\newcommand{\Cliff}{\mathop{\rm Cliff}\nolimits}

\def\ov{\overline}

\begin{document}
\maketitle

{\footnotesize{\bf Authors' address:} Departamento de Algebra, Universidad de Santiago
de Compostela. $15706$ Santiago de Compostela. Galicia. Spain. e-mail: {\tt
pedreira@zmat.usc.es}; \\ {\tt luisfg@usc.es}\\
{\bf Abstract:}  We study the projective normality of a linearly normal special scroll $R$ of degree
$d$ and speciality $i$ over a smooth curve $X$ of genus $g$. We relate it with the Clifford index of
the base curve $X$. If $d\geq 4g-2i-\Cliff(X)+1$, $i\geq 3$ and $R$ is smooth, we prove that the
projective normality of the scroll is equivalent to the projective normality of its directrix curve of
minimum degree.}

\vspace{0.1cm}

{\Large\bf Introduction.} 

Let $R\subset P^N$ be a linearly normal special scroll of genus $g$, speciality $i$ and degree $d$. We
know that $R$ has an associated ruled surface $\pi:S\lrw X$ with a linear system $|H|=|X_0+\b f|$, such
that $R$ is the image of $S$ by the map defined by $|H|$ (see \cite{fuentes}). We study the projective
normality of $R$, or equivalently, the normal generation of the invertible sheaf $\Te_S(H)$ on $S$.

In the previous paper \cite{normalidad1} we gave the bound $d\geq 4g-2i+1$ to reduce the problem of the
projective normality of the special scroll $R$ to the problem of the projective normality of the curve
of minimum degree.

In this paper we improve this result by using the Clifford index of the base curve $X$.  We will prove
the following theorem:

{\bf Theorem} {\em
Let  $R\subset \P^N$ be a smooth special linearly normal scroll of genus $g$, degree $d$ and speciality
$i\geq 3$. Let $X$ be the base curve of the scroll.  

If $d\geq 4g-2i-\Cliff(X)+1$, then:
\begin{enumerate}

\item $R$ has an unique special directrix curve $\ov{X_0}$. Moreover,
$\ov{X_0}$ is the curve of minimum degree, it is linearly normal and it has the speciality of $R$.

\item $R$ and $\ov{X_0}$ have the same speciality respect to hypersurfaces of degree $m$. In particular
the scroll is projectively normal if and only if the curve of minimum degree is projectively normal.

\end{enumerate} }

We will see that this result is optimal. We also study the particular case of $R$ being a cone. In
this case the bound $d\geq 4g-2i-\Cliff(X)+1$ is equivalent to the bound $d\geq 2g-2h^1(\L)-\Cliff(X)$,
that M. Green and R. Lazarsfeld gave in \cite{green} for a line bundle $\L$ on a curve.

We refer to \cite{fuentes} for a systematic development of the projective theory of scrolls and ruled
surfaces that we will use in this paper and to \cite{fuentes2} to study the special scrolls. In
the first section we recall some basic facts about the Clifford index of a curve, that we will use along
the paper.

\section{Preliminaries.}

A {\it geometrically ruled surface}, or simply a {\it
ruled surface}, will be a $\P^1$-bundle over a smooth curve $X$ of genus $g>0$. It will be
denoted by $\pi: S=\P(\E_0)\lrw X$. We will suppose that
$\E_0$ is a normalized sheaf and $X_0$ is the section of minimum self-intersection
that corresponds to the surjection $\E_0\lrw \Te_X(\e)\lrw 0$, $\bigwedge^2\E\cong
\Te_X(\e)$ and $e=-deg(\e)$ (see \cite{hartshorne},V, \S 2 and \cite{fuentes}).

If $|H|=|X_0+\b f|$ is a base-point-free linear system on a ruled surface $S$, $|H|$ defines a regular
map $\phi_H:S\lrw \P(H^o(\Te_S(H)^\vee))$. The image of $S$ is a scroll $R$. If $\phi_H$ is a birational
map we say that
$S$ and $H$ are the ruled surface and the linear system associated to the scroll $R$. We denote de image
of a curve $C\subset S$ by $\ov{C}\subset R$. The curve $\ov{X_0}$ is the curve of minimum degree of $R$.
It is embedded by the linear system $|\b+\e|$ on $X$.

Let $X$ be a smooth curve of genus $g\geq 2$. Let $\L$ be a line bundle on $X$. We define de Clifford
index of $\L$ by:
$$
\Cliff(\L)=deg(\L)-2(h^0(\L)-1)
$$
The Clifford index of the curve $X$ is defined by:
$$
\Cliff(X)=min\{\Cliff(\L)/h^0(\L)\geq 2, h^1(\L)\geq 2\}
$$
From this we have the following formula:

\begin{lemma}\label{rclifford}
If $\b$ is an effective special divisor such that $h^0(\Te_X(\b))\geq 2$ and $h^0(\Te_X(\b))\geq 2$ then
$$
h^0(\Te_X(\b))\leq \frac{deg(\b)-\Cliff(X)}{2}+1
$$ \qed
\end{lemma}

By the Clifford Theorem, $\Cliff(X)\geq 0$ with equality if and only if $X$ is hyperelliptic;
$\Cliff(X)=1$ if and only if either $X$ is trigonal or a smooth plane quintic. Furthermore, if $X$ is a
general curve of genus $g$ then $\gamma=[\frac{g-1}{2}]$, and in any event $\gamma\leq
[\frac{g-1}{2}]$.

Note that if $\b$ is a divisor such that $\Te_X(\b)$ provides the Clifford index of $X$, the linear
system $|\b|$ is base-point-free. In other case, if $P$ is a base point of $|\b|$, then $\b-P$ is a
divisor with a Clifford index less than $\Cliff(\Te_X(\b))$.

\section{Projective normality of a special scroll.}

\begin{prop}\label{directrizespecialminima2}
Let  $R\subset \P^N$ be a special linearly normal scroll of genus $g$, degree $d$ and speciality $i\geq
3$. Suppose that $R$ is not a cone. Let $X$ be the base curve of the scroll. 

If $d\geq 4g-2i-\Cliff(X)+1$, then $R$ has an unique special directrix curve $\ov{X_0}$.  Moreover,
$\ov{X_0}$ is the curve of minimum degree, it is linearly normal and it has the speciality of $R$.
\end{prop}
{\bf Proof:} Let  $S$ be the ruled surface and  $|H|=|X_0+\b f|$ the linear system corresponding to
the scroll $R$. Let $\gamma=\Cliff(X)$ be the Clifford index of $X$.

Since  $R$ is special, it has a special directrix curve (see \cite{fuentes2}) so the curve
$\ov{X_0}$ of minimum degree of the scroll  verifies $deg(\b+\e)\leq 2g-2$. Furthermore, we know that:
$$
d-2g+2+i=h^0(\Te_S(H))\leq h^0(\Te_X(\b+\e))+h^0(\Te_X(\b))\eqno(1)
$$
and
$$
i=h^1(\Te_S(H))\leq h^1(\Te_X(\b+\e))+h^1(\Te_X(\b))\eqno(2)
$$
Because $R$ is not a cone, $h^0(\Te_X(\b+\e))\geq 2$. We will prove that $deg(\b)\geq 2g+1$ and then we
will apply the Proposition $2.3$ of \cite{normalidad1}.

\begin{enumerate}

\item Suppose that $h^1(\Te_X(\b+\e))\geq 2$. Then we can apply the formula of Lemma 
\ref{rclifford} to the divisor $\b+\e$.

If $deg(\b)\leq 2g$, then we also can apply the Clifford Theorem  (\cite{hartshorne}, page 343) to
the divisor $\b$.  From the inequality $(1)$ we obtain:
$$
d-2g+2+i\leq \frac{deg(\b+\e)-\gamma}{2}+1+\frac{deg(\b)}{2}+1=\frac{d-\gamma}{2}+2
$$
and then  $d\leq 4g-2i-\gamma$ which contradicts the hypothesis.

\item Suppose that $h^1(\Te_X(\b+\e))\leq 1$. By hypothesis $i\geq 3$, so $h^1(\Te_X(\b))\geq 2$.

\begin{enumerate}

\item If $h^0(\Te_X(\b))\geq 2$ we can apply the formula of Lemma \ref{rclifford} to the
divisor $\b$ and the Clifford Theorem to the divisor $\b+\e$. From the inequality $(2)$ we obtain
$d\leq 4g-2i-\gamma$ which contradicts the hypothesis.

\item If $h^0(\Te_X(\b))\leq 1$, we have that:

$$
\begin{array}{rl}
{1\geq }&{h^0(\Te_X(\b))=deg(\b)-g+1+h^1(\Te_X(\b))}\\
{}&{\geq  deg(\b)-g+i}\\
\end{array} 
$$
Furthermore, by Nagata Theorem \cite{nagata} we know that $deg(\e)\leq g$, so:
$$
deg(\b+\e)\leq deg(b)+g\leq 2g-i+1 \eqno(3)
$$
On the other hand, from the inequality $(1)$ we have: 
$$
\begin{array}{rl}
{d-2g+2+i}&{\leq 1+deg(\b+\e)-g+1+h^1(\Te_X(\b+\e))\leq}\\
{}&{\leq 1+deg(\b+\e)-g+1+1}\\
\end{array}
$$
and because $d\geq 4g-2i-\gamma+1$,
$$
deg(\b+\e)\geq 3g-i-\gamma 
$$
Now, replacing the above expression at inequality $(3)$ we obtain:
$$
2g-i+1\geq 3g-i-\gamma \Longrightarrow \gamma\geq g-1
$$
but the Clifford index verifies $\gamma\leq [\frac{g-1}{2}]$. \qed

\end{enumerate}

\end{enumerate} 

\begin{rem}
{\em The inequality and the condition $i\geq 3$ are optimal in
the following way: 

Given a non
hyperelliptic smooth curve
$X$, let $\a$ be a divisor such that $\Te_X(\a)$ provides the Clifford index $\gamma$ of $X$. Let us
consider the ruled surface
$S=\P(\Te_X\oplus \Te_X(\a-\K))$. The linear system $|X_0+\k f|$ on $S$ is base-point-free and defines a
birrational map $\phi_H$. Let $R$ be the image of $S$ by the map $\phi_H$. The degree of $R$ is
$d=2g-2+deg(\a)$ and the speciality is $i=1+g-\frac{deg(a)+\gamma}{2}\geq 3$. From this:
$$
4g-2i-\gamma=4g-2-2g+deg(a)+\gamma-\gamma=d
$$
However, the scroll $R$ has two special directrix curves: $\ov{X_0}$ and $\ov{X_1}$ defined by the
linear systems $|\a|$ and $|\k|$ respectively.

On the other hand, we can also take the ruled surface $S=X\times P^1$ and the linear system $|X_0+\k
f|$. In this case, the corresponding scroll $R$ has degree $d=4g-4$ and speciality $2$. Since $X$ is non
hyperelliptic, $\gamma=\Cliff(X)\geq 1$ and $d\geq 4g-2i-\gamma$+1, but the scroll $R$ has a one
dimensional family of special directrix curves.}\qed
\end{rem}

\begin{prop}
Let  $R\subset \P^N$ be a special linearly normal scroll of genus $g$, degree $d$ and
speciality
$i\geq 3$. Suppose that $R$ \underline{is  a cone}. Let $X$ be the base curve of the scroll. 

If the unique singular point of $R$ is the vertex and $d\geq 4g-2i-\Cliff(X)+1$, then $R$ is
projectively normal.
\end{prop}
{\bf Proof:} We know that  $S=\P(\Te_X\oplus \Te_X(-\b))$ is the ruled surface
associated to $R$ and $R$ is given by the linear system $|X_1|=|X_0+\b f|$. Moreover, the degree of $R$
is $d=deg(\b)$ and the speciality is $i=g+h^1(\Te_X(\b))$ (see \cite{fuentes}).  

It is clear that $R$ is projectively normal iff $\Te_X(\b)$ is normally generated (see \cite{fuentes2}).
Since the unique singular point of $R$ is the vertex, the linear system $|\b|$ is very ample. Moreover,
$$
deg(\b)=d\geq 4g-2i-\Cliff(X)+1=2g-2h^1(\Te_X(\b))-\Cliff(X)+1
$$
Thus, we can apply the Green-Lazarsfeld Theorem (see \cite{green}) to the divisor $\b$ and we deduce that
the cone is projectively normal. \qed

\begin{rem}
{\em Note that the condition $d\geq 4g-2i-\Cliff(X)+1$ is optimal, because it is equivalent to the
inequality
$d\geq 2g-2h^1(\Te_X(\b))-\Cliff(X)+1$ for the hyperplane section of the cone.  This condition is
the best possible for the projective normality of line bundles on curves (see \cite{green}).} \qed
\end{rem}

\begin{teo}\label{principal}
Let  $R\subset \P^N$ be a smooth special linearly normal scroll of genus $g$, degree $d$ and speciality
$i\geq 3$. Let $X$ be the base curve of the scroll.  

If $d\geq 4g-2i-\Cliff(X)+1$, then:
\begin{enumerate}

\item $R$ has an unique special directrix curve $\ov{X_0}$. Moreover,
$\ov{X_0}$ is the curve of minimum degree, it is linearly normal and it has the speciality of $R$.

\item $R$ and $\ov{X_0}$ have the same speciality respect to hypersurfaces of degree $m$. In particular
the scroll is projectively normal if and only if the curve of minimum degree is projectively normal.

\end{enumerate}

\end{teo}
{\bf Proof:} Let  $S$ be the ruled surface and  $|H|=|X_0+\b f|$ the linear system corresponding to
the scroll $R$. Let $\gamma=\Cliff(X)$ be the Clifford index of $X$.

The first assertion is the Proposition \ref{directrizespecialminima2}. From the proof of this Proposition
we also know that $deg(\b)\geq 2g+1$. 

To prove the second statement we will apply the Proposition $2.1$ of \cite{normalidad1}. We will see that:
$$
s(\b+\e,\stackrel{i}{\ldots},\b+\e,\b,\stackrel{k-i}{\ldots},\b)=0 \mbox{ for all $i$, with $0\leq i\leq
k-1$}
$$
Reasoning as in the proof of the Theorem $2.4$ of \cite{normalidad1} it is sufficient to see that
$s(\b,\b+\e)=0$ and in particular we only have to prove that (see Lemma $1.5$ in \cite{normalidad1}):
$$
h^1(\Te_X(\b-(\b+\e)))< h^0(\Te_X(\b+\e))-1
$$
We distinguish two cases:

\begin{enumerate}

\item Suppose that $h^1(\Te_X(-\e))\leq 1$. It is sufficient to prove that $h^0(\Te_X(\b+\e))\geq 3$.
But this follows from the smoothness of the scroll: 

\begin{enumerate}

\item If $h^0(\Te_X(\b+\e))=0$ the scroll is a cone. 

\item $h^0(\Te_X(\b+\e))=1$ can not occur, because $\b+\e$ is base-point-free.

\item If $h^0(\Te_X(\b+\e))=2$ the directrix curve of minimum degree is a line. Since the scroll is not
rational, it must be a singular curve of the scroll.

\end{enumerate}

\item Suppose that $h^1(\Te_X(-\e))\geq 2$:

\begin{enumerate}

\item If $h^0(\Te_X(-\e))\leq 1$, then by Riemann-Roch
Theorem $h^1(\Te_X(-\e))\leq g-e$.  Moreover, we know that
$deg(\b)\geq 2g+1$ and from this:
$$
h^0(\Te_X(\b+\e))-1=b-e-g+i\geq 2g+1-e-g+i> g-e\geq h^1(\Te_X(-\e))
$$

\item If $h^0(\Te_X(-\e))\geq 2$ we can apply the formula of Lemma \ref{rclifford} to the divisor
$-\e$:
$$
h^0(\Te_X(-\e))\leq \frac{e-\gamma}{2}-1\mbox{, or equivalently, } h^1(\Te_X(-\e))\leq
g-\frac{e+\gamma}{2}
$$
Furthermore,
$$
h^0(\Te_X(\b+\e))-1=deg(\b)-e-g+1+i-1=deg(\b)-e-g+i
$$
By hypothesis $d\geq 4g-2i-\gamma+1$, so $2deg(\b)-e\geq 4g-2i-\gamma+1$ and $deg(\b)\geq
\frac{4g-2i-\gamma+1+e}{2}$. Then:
$$
\begin{array}{rl}
{h^0(\Te_X(\b+\e))-1}&{\geq \frac{4g-2i-\gamma+1+e}{2}-e-g+i\geq g-\frac{e+\gamma}{2}+\frac{1}{2}\geq}\\
{}&{}\\
{}&{\geq h^1(\Te_X(-\e))+\frac{1}{2}>h^1(\Te_X(-\e))}\\
\end{array}
$$

\end{enumerate}

\end{enumerate}

\qed

\end{document}